\documentclass{article}

       \newcommand{\R}{\ensuremath{\sf R\hspace*{-1.0ex}\rule{0.15ex}%
       {1.5ex}\hspace*{0.9ex}}}

       \newcommand{\LC}{{\sf \hspace*{2ex}\rule{0.15ex}%
       {1.5ex}\hspace*{-1.5ex}\rule{1.5ex}{0.15ex}\hspace*{0.5ex}}}
       \newcommand{\RC}{{\sf \hspace*{0.7ex}\rule{0.15ex}%
       {1.5ex}\hspace*{-0.15ex}\rule{1.5ex}{0.15ex}\hspace*{0.5ex}}}
       \newcommand{\be}{\begin{equation}}
       \newcommand{\ee}{\end{equation}}
       \newcommand{\bv}[1]{\mathrm{\mathbf{#1}}}

\usepackage[dvips]{graphicx}
\usepackage{verbatim}
\begin{document}
\title{Euclidean Geometric Objects in the Clifford Geometric Algebra of \{Origin, 3-Space, Infinity\}} 
\author{Eckhard M.S. Hitzer}
\date{1 Dec. 2003 (corrections: 14 July 2004)}
\maketitle

\begin{abstract}
This paper concentrates on the homogeneous (conformal) model of
Euclidean space (\textit{Horosphere}) with subspaces that intuitively 
correspond to Euclidean geometric objects in three dimensions. 
Mathematical details of the construction and (useful) parametrizations of the 3D
Euclidean object models are explicitly demonstrated in order to show how 
3D Euclidean information on positions, orientations and radii can be extracted. 
\end{abstract}

\section{Introduction}

The Clifford geometric algebra of three dimensional (3D) Euclidean space with vectors
\be
 \bv{p}=p_1\bv{e}_1+p_2\bv{e}_2+p_3\bv{e}_3.
\ee
given in terms of an orthonormal basis $\{\mathbf{e}_1, \mathbf{e}_2, \mathbf{e}_3  \}$, 
nicely encodes the algebra of 3D subspaces with algebraic basis
\be
\{1, \mathbf{e}_1, \mathbf{e}_2, \mathbf{e}_3, 
    \bv{i}_1=\mathbf{e}_2\mathbf{e}_3, 
    \bv{i}_2=\mathbf{e}_3\mathbf{e}_1, 
    \bv{i}_3=\mathbf{e}_1\mathbf{e}_2,
    i = \mathbf{e}_1 \mathbf{e}_2 \mathbf{e}_3 \}, 
\ee
providing geometric multivector product expressions of 
rotations and set theoretic opeations.\cite{DH:NFCM} 
But in this framework line and plane subspaces always contain the origin. 

The homogeneous (conformal) model of 3D Euclidean space in the Clifford geometric algebra
$\R_{4,1}$ provides a way out.\cite{{GS:Wachter}, {HLR:GS1}} Here positions of points, lines and planes, etc. off the 
origin can be naturally encoded. 
Other advantages are the unified treatment of rotations and translations and ways to 
encode point pairs, circles and spheres. The creation of such elementary geometric objects
simply occurs by algebraically joining a minimal number of points in the object subspace.
The resulting multivector expressions completely encode in their components
positions, orientations and radii. 

The geometric algebra $\R_{4,1}$ can be intuitively
pictured as the algebra of origin $\bar{\mathbf{n}}$, Euclidean 3D space and 
infinity $\mathbf{n}$, where origin and infinity are represented by additional 
linearly independent null-vectors. 
\be
  \{\bar{\mathbf{n}}, \mathbf{e}_1, \mathbf{e}_2, \mathbf{e}_3, \mathbf{n} \},
  \,\,\, \mathbf{n}^2=\bar{\mathbf{n}}^2=0.
  \label{eq:5Db}
\ee
This algebra
seems most suitable for applications in computer graphics, robotics and other 
fields.\cite{{DLL:CGeo}, GS:GeoComp}

This paper concentrates on giving explicit details for the construction of fundamental 
geometric objects in this model, detailing how the 3D geometric information can be
extracted. How the simple multivector representations of these objects 
can be manipulated in order to move them in three dimensions and to express set theoretic
operations of union (join), intersection (meet), projections and rejections is described
in~\cite{{EH:KWA}, {HLR:GS1}, {DLL:CGeo}}.

This algebraic encoding of geometric objects and their manipulations strongly suggests
an object oriented software implementation. This does allow computers to calculate
with this algebra and provides programmers with the means to most suitably represent 
fundamental geometric objects, their 3D properties and ways (methods) to manipulate 
these objects. This happens on a higher algebraic level, so that the programmer actually
is freed of the need to first investigate suitable less intuitive matrix 
representations.\cite{{EH:KWA}, {EH:GAP}, {CP:CLU}, {LD:gaigen}, {LD:GAviewer}}

The mathematical notation used in this paper for multivectors, geometric products and 
products derived from the geometric product is fairly standard: Italic capital letters
are used for conformal vectors and multivectors, bold lower case vectors represent 
Euclidean vectors (except the origin and infinity null vectors), the wedge $\wedge$ 
signifies the outer product, the asterisk $\ast$ the scalar product, the angles $\LC$
and $\RC$ left and right contractions and mere juxtaposition the full geometric product
of multivectors. Other notations are explicitly defined where they are 
used.\cite{LD:IPofGA, EH:GPderP}

\section{3D Information in Homogeneous Objects}

We will see how homogenous multivectors completely encode positions, 
directions, moments and radii of the corresponding three dimensional (3D) objects
in Euclidean space. An overview of this is give in Table \ref{tb:3Dinf}.
\begin{table}
\begin{center}
\begin{tabular}{c|c}
  homogeneous object   & 3D information  
  \\ \hline \hline
  point $X$ & position $\mathbf{x}$               
  \\ \hline
  point pair $P_1\wedge P_2$ & positions $\mathbf{p}_1,\mathbf{p}_2$ 
  \\ \hline
  line & direction vector, moment bivector
  \\ \hline
  cirlce & plane bivector, center, radius
  \\ \hline
  plane & plane bivector, location vector
  \\ \hline
  sphere & center, radius
\end{tabular}
\caption{3D geometric information in homogeneous objects.
The left column lists the homogeneous multivectors, 
that represent the geometric objects.}
\label{tb:3Dinf}
\end{center}
\end{table}
In the rest of this section we will look at the details of extracting the encoded
3D information from each homogeneous multivector object. Where suitable, we will also
give useful alternative parametrizations of homogeneous multivector 
objects.\footnote{Dorst and Fontijne\cite{DF:ITM03proc} give similar parametrizations,
but with the rather strong simplification, that objects are centered at the origin 
$\bar{\mathbf{n}}$ or contain the origin $\bar{\mathbf{n}}$, e.g. $C=\bar{\mathbf{n}}$
in eq. (\ref{eq:PPcmid}), etc. The general and explicit 
formulas presented in the following, seem to appear nowhere else in the published 
literature so far.}

\subsection{Point and Pair of Points}

The Euclidean position $\mathbf{p}$ of a conformal point 
\be
  P=\bv{p}+\frac{1}{2}p^2\bv{n}+\bv{\bar{n}}
  \label{eq:confP}
\ee
 is obtained with 
the help of the (additive\cite{HLR:GS1}) \textit{conformal split}, which is an example of a 
rejection\cite{DH:CAGC}:
The conformal point vector $P$ is rejected off the Minkowski plane represented by the
bivector $N=\bv{n} \wedge \bv{\bar{n}}$
\be
  \mathbf{p} = (P\wedge N)N.
\ee
Equation (\ref{eq:confP}) shows how to achieve the opposite, i.e. how to
get back to the conformal point $P$ from just knowing the Euclidean position $\mathbf{p}$.

The Euclidean positions $\mathbf{p}_1,\mathbf{p}_2$ of a pair of points represented
by the conformal bivector 
\begin{eqnarray}
  V_2 & = & P_1\wedge P_2
  \nonumber \\
      & = & \mathbf{p}_1\wedge\mathbf{p}_2
%
                     +\frac{1}{2}(p_2^2\mathbf{p}_1 - p_1^2\mathbf{p}_2)\mathbf{n}
                     -(\mathbf{p}_2-\mathbf{p}_1)\bar{\mathbf{n}}
                     +\frac{1}{2}(p_1^2- p_2^2)N
  \nonumber \\
      & = & \bv{b} -\frac{1}{2}\mathbf{v}\mathbf{n}     
                      -\mathbf{u}\bar{\mathbf{n}}
                      -\frac{1}{2}\gamma N
\end{eqnarray}
can be fully reconstructed from the components of $V_2$. 
%
%
We assume without restricting the generality, 
that $p_1 = \sqrt{\mathbf{p}_1^2}\geq p_2 = \sqrt{\mathbf{p}_2^2}$.
Given any conformal bivector $V_2$ with
components $\bv{b}$ (a Euclidean bivector), $\mathbf{u}$ and  $\mathbf{v}$ (Euclidean vectors of 
length $u=\sqrt{\mathbf{u}^2}$ and $v=\sqrt{\mathbf{v}^2}$), and $\gamma$ (a real scalar), 
the calculation works as follows
\begin{eqnarray}
  && \sigma = \frac{1}{2}\gamma^2-\mathbf{u} \ast \mathbf{v}, \hspace{5mm}
     \rho = \sqrt{\sigma^2-u^2v^2}, 
  \label{eq:pair1} \\
  && p_1= \frac{\sqrt{\sigma +\rho}}{u}, \hspace{5mm} 
     p_2= \frac{\sqrt{\sigma -\rho}}{u}, \\
  && \mathbf{p}_1 = p_1\frac{p_1^2 \mathbf{u}+\mathbf{v}}{\left| p_1^2 \mathbf{u}+\mathbf{v} \right|},
     \hspace{5mm}
%
%
%
     \mathbf{p}_2 = p_2\frac{p_2^2 \mathbf{u}+\mathbf{v}}{\left| p_2^2 \mathbf{u}+\mathbf{v} \right|}. 
  \label{eq:pair3}
\end{eqnarray}
This calculation is the full solution (of two conformal points $X=P_1,P_2$) to the equation
\be
  V_2 \wedge X = 0, X^2=0.
\ee

We can further view conformal point pairs as one-dimensional circles and arrive
thereby at another highly useful 
characterization\footnote{This characterization is e.g. very useful for investigating
the full (real and virtual) meet of two circles, or of a straight line and a 
circle.\cite{EH:GIAEpres}}: 
\be
   P_1\wedge P_2  
   = 2r \{\bv{\hat{p}}\wedge \bv{c} 
       + \frac{1}{2}[(c^2 + r^2)\bv{\hat{p}} -2\bv{c}\ast\bv{\hat{p}}\,\,\bv{c} ]\mathbf{n}
       + \bv{\hat{p}}\bv{\bar{n}}+\bv{c}\ast\bv{\hat{p}}N \},
   \label{eq:pp1Dc}
\ee
with the "radius" $r$ defined as half the Euclidean point pair distance, 
$\bv{\hat{p}}$ a unit vector pointing from $\mathbf{p}_2$ to $\mathbf{p}_1$,
and $\bv{c}$ the Euclidean midpoint (center) of the point pair:

\be 
     2r =\, \mid\mathbf{p}_1-\mathbf{p}_2 \mid,     \,\,\, 
     \bv{\hat{p}}=\frac{\mathbf{p}_1-\mathbf{p}_2}{2r},  \,\,\,
     \bv{c}=\frac{\mathbf{p}_1+\mathbf{p}_2}{2}.
     \label{eq:rpc}
\ee
In case that the straight line defined by the point pair contains the origin,
i.e. for $\bv{\hat{p}}\wedge \bv{c} =0 \,\, (\bv{\hat{p}}\, \| \, \bv{c})$ we get the simplified form
\be
 P_1\wedge P_2=2r \{C-r^2 \bv{n}\} \bv{\hat{p}} N.
\ee
In case that the Euclidean midpoint vector $\bv{c}$ is perpendicular to $\bv{\hat{p}}$ 
$(\bv{\hat{p}}\perp \bv{c})$, i.e. if $\bv{\hat{p}}\ast \bv{c}=0$ we get
\be
  P_1\wedge P_2=-2r \{C+r^2 \bv{n}\}\bv{\hat{p}}.
\ee
In both cases we used the conformal representation of the midpoint as
\be
 C=\bv{c}+\frac{1}{2}c^2\bv{n}+\bv{\bar{n}}.
 \label{eq:PPcmid}
\ee

\subsection{Lines}

Given two conformal points $P_1$ and $P_2$ the conformal trivector
\be
  V_{line} = P_1\wedge P_2 \wedge \mathbf{n}
           = \mathbf{p}_1\wedge\mathbf{p}_2 \wedge \mathbf{n}
             + (\mathbf{p}_2-\mathbf{p}_1) N
           = \bv{m}\, \mathbf{n} + \mathbf{d} N
\ee
conveniently consists of the defining entities of the Euclidean line 
through $\mathbf{p}_1$ and $\mathbf{p}_2$. 
The Euclidean bivector $\bv{m}$ 
represents the moment and the Euclidean vector $\mathbf{d}$ the direction of the line.
Using $\bv{m}$ and $\mathbf{d}$ we can give the parametric form of the line as
\be
  \mathbf{x} = (\bv{m}+\alpha)\mathbf{d}^{-1}, \alpha \in \R.
  \label{eq:Lsol}
\ee
All points $X=\mathbf{x}+\frac{1}{2}x^2\mathbf{n}+ \bar{\mathbf{n}}$ with the $\mathbf{x}$
as specified in (\ref{eq:Lsol}) represent the full solution to the problem
\be
  V_{line}\wedge X =0, X^2=0.
\ee

The one-dimensional circle representation of point pairs (\ref{eq:pp1Dc}) immediately leads
to a second often useful parametrization of lines as
\be
 P_1 \wedge P_2 \wedge \bv{n} 
   = 2r \bv{\hat{p}}\wedge C \wedge n
   = 2r \{ \bv{\hat{p}}\wedge \bv{c}\,\,\bv{n} - \bv{\hat{p}}N \}.
   \label{eq:line1Dc}
\ee
It is important to note that the conformal point $C$ in eq. (\ref{eq:line1Dc}) 
does not need to be the midpoint of the point pair. Any conformal point on the 
straight line $P_1 \wedge P_2 \wedge \bv{n} $ can take the place of $C$
in eq. (\ref{eq:line1Dc}). $\bv{\hat{p}}$ and $r$ are defined as in eq. (\ref{eq:rpc}).

\subsection{Circles}

\label{ssc:circle}
General conformal trivectors of the form
\be
  V_3 = P_1\wedge P_2 \wedge P_3,
  \label{eq:V3P3}
\ee
with conformal points $P_1, P_2$ and $P_3$ represent Euclidean circles through the
corresponding Euclidean points $\mathbf{p}_1, \mathbf{p}_2$ and $\mathbf{p}_3$. The
equation for all points $X$ on such a circle is again given as
\be
  V_3\wedge X =0, X^2=0.
\ee
In order to clearly interpret and apply $V_3$ and its various components we will explicitly 
insert the three points
\be
  P_1 = \mathbf{p}_1 + \frac{1}{2}p_1^2 \mathbf{n} + \bar{\mathbf{n}},
  P_2 = \mathbf{p}_2 + \frac{1}{2}p_2^2 \mathbf{n} + \bar{\mathbf{n}},
  P_3 = \mathbf{p}_3 + \frac{1}{2}p_3^2 \mathbf{n} + \bar{\mathbf{n}}.
  \label{eq:cpoints}
\ee
The conformal circle trivector becomes
\begin{eqnarray}
  V_3 &=&  \mathbf{p}_1\wedge \mathbf{p}_2 \wedge \mathbf{p}_3 
  \nonumber \\
  && +\frac{1}{2}( p_1^2\mathbf{p}_2 \wedge \mathbf{p}_3
               +p_2^2\mathbf{p}_3 \wedge \mathbf{p}_1
               +p_3^2\mathbf{p}_1 \wedge \mathbf{p}_2) \mathbf{n}
  \nonumber \\
  &&+(\mathbf{p}_2 \wedge \mathbf{p}_3+\mathbf{p}_3 \wedge \mathbf{p}_1
      +\mathbf{p}_1 \wedge \mathbf{p}_2) \bar{\mathbf{n}}
  \nonumber \\
  &&
  +\frac{1}{2}\{\mathbf{p}_1(p_2^2-p_3^3) +\mathbf{p}_2(p_3^2-p_1^3)+\mathbf{p}_3(p_1^2-p_2^3)\}N
  \label{eq:V3p}
\end{eqnarray}
The Euclidean bivector component factor of $\bar{\mathbf{n}}$ 
\begin{eqnarray}
  I_c & = & -\{[V_3+(V_3\ast i)i] \wedge \mathbf{n}\} N 
      \nonumber \\
      & = & \mathbf{p}_2 \wedge \mathbf{p}_3+\mathbf{p}_3 \wedge \mathbf{p}_1
            +\mathbf{p}_1 \wedge \mathbf{p}_2 
      \nonumber \\
      & = & (\mathbf{p}_1 - \mathbf{p}_2)\wedge(\mathbf{p}_2-\mathbf{p}_3)
  \label{eq:Ic}
\end{eqnarray}
is obviously parallel to the plane (of the Euclidean circle) through
$\mathbf{p}_1, \mathbf{p}_2$ and $\mathbf{p}_3$. 
Assuming the Euclidean center vector of the circle to be $\mathbf{c}$ and the radius $r$, we can
rewrite (\ref{eq:cpoints}) as
\be
  P_k = \mathbf{c} + r \mathbf{r}_k 
        + \frac{1}{2}(c^2 + r^2 +2r \mathbf{c}\ast \mathbf{r}_k)\mathbf{n} + \bar{\mathbf{n}},\,\,\,
        \mathbf{r}_k^2=1,\,\,\,
        k = 1,2,3.
\ee
The three vectors  $\mathbf{r}_k$ are unit length vectors pointing from the circle center $\mathbf{c}$
to the three points $\mathbf{p}_1, \mathbf{p}_2$ and $\mathbf{p}_3$, respectively. Replacing the
 $P_k$ in (\ref{eq:V3p}) accordingly we get after doing some algebra the simplified form
\be
  V_3 = \mathbf{c}\wedge I_c 
        + [\frac{1}{2}(r^2+c^2)I_c - \mathbf{c}(\mathbf{c}\LC I_c)]\mathbf{n}
        + I_c \bar{\mathbf{n}} - (\mathbf{c}\LC I_c) N.
  \label{eq:V3cr}
\ee
We see that the three vectors $\mathbf{r}_k, k=1,2,3$ do no longer occur explicitly. They
enter equation (\ref{eq:V3cr}) only by defining the orientation of the circle plane $I_c$ in 
(\ref{eq:Ic}). 

If we assume only to know $V_3$ as outer product (\ref{eq:V3P3}) of three general conformal points
we can now extract the radius $r$ by calculating
\be
  r^2 = -\frac{V_3^2}{I_c^2}.
\ee
We can decompose the center vector $\mathbf{c}$ by way of projection and
rejection into components parallel and perpendicular to the circle plane
\begin{eqnarray}
  \mathbf{c} &=& \mathbf{c}_{\parallel}+\mathbf{c}_{\perp}, \\
  \mathbf{c}_{\parallel}&=& (\mathbf{c}\LC I_c^{-1}) I_c
  = (\mathbf{c}\LC I_c) I_c^{-1} 
  \stackrel{(\ref{eq:V3cr})}{=}-[(V_3 \RC {\mathbf{n}})\RC \bar{\mathbf{n}}] I_c^{-1}, 
  \label{eq:cpar}  \\
  \mathbf{c}_{\perp}&=& (\mathbf{c}\wedge I_c^{-1})I_c
  = (\mathbf{c}\wedge I_c) I_c^{-1} 
  \stackrel{(\ref{eq:V3cr})}{=} -(V_3\ast i)i I_c^{-1}.
  \label{eq:cperp}
\end{eqnarray}
The Euclidean circle center vector can hence be extracted from any $V_3$ as
\be
  \mathbf{c} = \mathbf{c}_{\parallel}+\mathbf{c}_{\perp}
    \stackrel{(\ref{eq:cpar}),(\ref{eq:cperp}),(\ref{eq:Ic})}{=} 
    -\frac{(V_3 \RC {\mathbf{n}})\RC \bar{\mathbf{n}} + (V_3\ast i)i}{I_c}
\ee

Inserting the decomposition $\mathbf{c} = \mathbf{c}_{\parallel}+\mathbf{c}_{\perp}$ 
we get the following expression for the circle trivector
\begin{eqnarray}
  V_3 &=& \mathbf{c}_{\perp} I_c 
        + [\frac{1}{2}(r^2-c^2)I_c + \mathbf{c}\mathbf{c}_{\perp} I_c]\mathbf{n}
        + I_c \bar{\mathbf{n}} - \mathbf{c}_{\parallel} I_c N
  \nonumber \\
      &=& \{\mathbf{c}_{\perp} N
        + [\frac{1}{2}(r^2-c^2) + \mathbf{c}\mathbf{c}_{\perp} ]\mathbf{n}
        - \bar{\mathbf{n}} - \mathbf{c}_{\parallel} \}I_c N
   \nonumber \\
      &=& \{ - \mathbf{c}_{\parallel} - \frac{1}{2}c_{\parallel}^2\mathbf{n} - \bar{\mathbf{n}}
           + \frac{1}{2}r^2\mathbf{n} 
           + \mathbf{c}_{\perp} N
           + [-\frac{1}{2}c_{\perp}^2 + \mathbf{c}\mathbf{c}_{\perp}]\mathbf{n} \}I_c N
\end{eqnarray}
In the case that the circle plane includes the origin ($\mathbf{c}_{\perp}=0$) we are left 
with
\be
  V_3 = -[C - \frac{1}{2}r^2\mathbf{n}] I_c N,
\ee
and can extract the conformal center
\be
  C=\mathbf{c} + \frac{1}{2}c^2\mathbf{n} + \bar{\mathbf{n}}
\ee
 simply as
\be
  C=-\frac{V_3}{I_c N}+\frac{1}{2}r^2\mathbf{n}.
\ee

\subsection{Planes}

Given three conformal points $P_1, P_2$ and $P_3$ as in (\ref{eq:cpoints}) 
the conformal 4-vector
\begin{eqnarray}
  V_{plane}&=&P_1\wedge P_2\wedge P_3 \wedge \mathbf{n} \nonumber \\
           &=& \mathbf{p}_1\wedge \mathbf{p}_2 \wedge \mathbf{p}_3 \wedge \mathbf{n} 
           \nonumber\\
           &&-(\mathbf{p}_2 \wedge \mathbf{p}_3+\mathbf{p}_3 \wedge \mathbf{p}_1
                +\mathbf{p}_1 \wedge \mathbf{p}_2) N
\end{eqnarray}
represents the plane through the Euclidean points 
$\mathbf{p}_1, \mathbf{p}_2$ and $\mathbf{p}_3$.
The Euclidean bivector component factor of $N$
\be
%
%
%
  I_p = -(V_{plane} \mathbf{n})\RC \bar{\mathbf{n}}
      = \mathbf{p}_2 \wedge \mathbf{p}_3+\mathbf{p}_3 \wedge \mathbf{p}_1
                +\mathbf{p}_1 \wedge \mathbf{p}_2
      = (\mathbf{p}_1 - \mathbf{p}_2)\wedge(\mathbf{p}_2-\mathbf{p}_3)
\ee
gives the orientation of the plane in the Euclidean space. 
This allows us to rewrite $V_{plane}$ as
\be
  V_{plane} = (\mathbf{d} \wedge I_p) \mathbf{n} - I_p N
            = \mathbf{d} I_p \mathbf{n} - I_p N,
\ee
where $\mathbf{d}$ represents the Euclidean distance vector from the origin to the
plane, itself perpendicular to the plane. The Euclidean distance vector can be extracted from
$V_{plane}$ by
\be
  \mathbf{d}= (V_{plane} \wedge \bar{\mathbf{n}})I_p^{-1}N.
\ee
The equation for all points $X$ on the plane is again given as
\be
  V_{plane}\wedge X =0, X^2=0.
\ee

Somewhat in analogy of the relation of point pairs as one-dimensional circles (\ref{eq:pp1Dc}) 
and the resulting alternative parametrization of lines (\ref{eq:line1Dc}), an alternative
parametrization of planes by means of a general conformal point $C$ on the plane is
possible
\be
  P_1 \wedge P_2 \wedge P_3 \wedge \bv{n}
    = C\wedge I_c \wedge \bv{n}
    = \bv{c}\wedge I_c\bv{n} - I_c N
\ee
For $\bv{c}\wedge I_c=0$ (origin $\bv{\bar{n}}$ in plane) we get
\be
  P_1 \wedge P_2 \wedge P_3 \wedge \bv{n}
    = -I_c N.
\ee

\subsection{Spheres}

\label{ssc:sphere}
General conformal 4-vectors of the form 
\be
  V_4 = P_1\wedge P_2\wedge P_3 \wedge P_4
\ee
with conformal points 
\be 
  P_k= \mathbf{p}_k + \frac{1}{2}p_k^2 \mathbf{n} + \bar{\mathbf{n}}, k=1,2,3,4
  \label{eq:cp4}
\ee
represent Euclidean spheres through the corresponding Euclidean points $\mathbf{p}_k, k=1,2,3,4$.
The equation for all points $X$ on the sphere is again given as
\be
  V_4\wedge X =0, X^2=0.
\ee
Inserting (\ref{eq:cp4}) explicitly in $V_4$ yields
\begin{eqnarray}
  V_4 &=& -\frac{1}{2}( p_1^2\mathbf{p}_{234}
                       +p_2^2\mathbf{p}_{314} 
                       +p_3^2\mathbf{p}_{124} 
                       +p_4^2\mathbf{p}_{132} ) \mathbf{n}
  \nonumber \\
      && -(             \mathbf{p}_{234}
                       +\mathbf{p}_{314} 
                       +\mathbf{p}_{124} 
                       +\mathbf{p}_{132})\bar{\mathbf{n}}
  \nonumber \\
  &&
  +\frac{1}{2}\{(p_2^2-p_3^2)\mathbf{p}_{14} 
               +(p_3^2-p_1^2)\mathbf{p}_{24}
               +(p_1^2-p_2^2)\mathbf{p}_{34}
  \nonumber \\
  &&
  \,\,\,\,\,\,      +(p_1^2-p_4^2)\mathbf{p}_{23}
               +(p_2^2-p_4^2)\mathbf{p}_{31}
               +(p_3^2-p_4^2)\mathbf{p}_{12}  \}N,
  \label{eq:V4p}
\end{eqnarray}
with the abbreviations
\be
  \mathbf{p}_{kl} =\mathbf{p}_k\wedge \mathbf{p}_l,  \,\,\,
  \mathbf{p}_{klm}=\mathbf{p}_k\wedge \mathbf{p}_l\wedge \mathbf{p}_m, \,\,\,
  k,l,m \in \{1,2,3,4\}.
\ee
The $\bar{\mathbf{n}}$ factor component 
\begin{eqnarray}
%
%
%
%
  i_s &=& -(V_4 \wedge \mathbf{n}) N 
  \nonumber \\
      &=& -(\mathbf{p}_{234} +\mathbf{p}_{314} +\mathbf{p}_{124} +\mathbf{p}_{132})
  \nonumber \\
      &=& (\mathbf{p}_1-\mathbf{p}_2)
          \wedge(\mathbf{p}_2-\mathbf{p}_3)
          \wedge(\mathbf{p}_3-\mathbf{p}_4)
\end{eqnarray}
is a Euclidean pseudoscalar, i.e. proportional to $i$.
Similar to the discussion of the circle, assuming the Euclidean center vector of the sphere 
to be $\mathbf{c}$ and the radius $r$, we can
rewrite (\ref{eq:cp4}) as
\be
  P_k = \mathbf{c} + r \mathbf{r}_k 
        + \frac{1}{2}(c^2 + r^2 +2r \mathbf{c}\ast \mathbf{r}_k)\mathbf{n} + \bar{\mathbf{n}},\,\,\,
        \mathbf{r}_k^2=1,\,\,\,
        k = 1,2,3,4.
\ee
Replacing the $P_k, k=1,2,3,4$ in (\ref{eq:V4p}) accordingly we get after doing lots of algebra
\begin{eqnarray}
  V_4 &=& \frac{1}{2}(r^2-c^2)i_s \mathbf{n}+i_s \bar{\mathbf{n}} + \mathbf{c}i_s N
  \nonumber \\
      &=& (\mathbf{c}+\frac{1}{2}c^2\mathbf{n}+ \bar{\mathbf{n}} - \frac{1}{2}r^2\mathbf{n}) i_s N
  \nonumber \\
      &=& (C- \frac{1}{2}r^2\mathbf{n}) i_s N,
  \label{eq:V4conf}
\end{eqnarray}
where $C$ represents the conformal center of the sphere. An important relationship used in the 
derivation of (\ref{eq:V4conf}) is 
\be
  i_s = (\mathbf{p}_1-\mathbf{p}_2)
          \wedge(\mathbf{p}_2-\mathbf{p}_3)
          \wedge(\mathbf{p}_3-\mathbf{p}_4)
       = r^3 (\mathbf{r}_1-\mathbf{r}_2)
          \wedge(\mathbf{r}_2-\mathbf{r}_3)
          \wedge(\mathbf{r}_3-\mathbf{r}_4).
\ee

The elegant form (\ref{eq:V4conf}) of $V_4$ makes it easy to extract the radius and the center from
any general conformal (sphere) 4-vector:
\be
  r^2 = \frac{V_4^2}{(V_4\wedge \mathbf{n})^2}, \hspace{5mm}
    C = \frac{1}{2}r^2\mathbf{n}+\frac{V_4}{-V_4\wedge \mathbf{n}}.
\ee

\section{Conclusions}

We explained how to algebraically construct conformal (homogeneous)  subspaces with
very intuitive Euclidean interpretations.\footnote{They are e.g. implemented, together with 
algebraic expressions for arbitrary translations and rotations, and for subspace operations 
of union (join), intersection (meet), projection and rejection as methods in the 
GeometricAlgebra Java package.\cite{{EH:KWA}, {EH:GAP}, {HU:HMJI}} }

We then analysed in detail how the joining of conformal points yields explicit expressions
for points, pairs of points, lines, circles, planes ans spheres.
After that we showed how the Euclidean 3D information of positions, orientations and radii, 
etc. can be 
extracted.\footnote{These formulas presicely yield the optimal mathematical 
structure of the related Java methods each geometric object is to have e.g. 
in the GeometricAlgebra Java package implementation.~\cite{EH:GAP}} 
In some cases useful alternative parametrizations were given.
Applications of these alternative parametrizations can e.g. be found in~\cite{EH:GIAEpres}.

\section*{Acknowledgement}

Soli Deo Gloria. The author does want to thank his wife and children. He further thanks
S. Krausshaar and P. Leopardi for organizing the ICIAM 2003 Clifford mini-symposium.

\end{document}